\documentclass[12pt]{article}

\usepackage{amssymb,amsthm,amsmath,latexsym}
\usepackage{url}
\usepackage{lscape}
\newtheorem{thm}{Theorem}[section]

\newtheorem{cor}[thm]{Corollary}

\theoremstyle{remark}

\title{The existence of partitioned balanced tournament designs}

\author{%
Makoto Araya\thanks{Department of Computer Science, Shizuoka University}\quad and \quad Naoya Tokihisa\thanks{
Connected Technology Planning Group, Mobility Service Department, New Mobility Service, SUZUKI MOTOR CORPORATION. The work of the second author was carried out at
Graduate School of Integrated Science and Technology, Shizuoka University.}
}

\begin{document}

\maketitle

\begin{abstract}
E. R. Lamken prove in \cite{pbtd} that 
there exists a partitioned balanced tournament design of side $n$, PBTD($n$), 
for $n$ a positive integer, $n \ge 5$, except possibly for $n \in \{9,11,15\}$.
In this article, we show the existence of PBTD($n$) for $n \in \{9,11,15\}$.
%It is not difficult to prove that there is no PBTD($n$) for $n \le 4$.
As a consequence, the existence of PBTD($n$) has been completely determined.

\end{abstract}

\section{Introduction}

%E. R. Lamken prove in \cite{pbtd} that the existence of the PBTD($n$) for $n$ a positive integer, $n \ge 5$, except possibly for $n \in \{9,11,15\}$.
%In this article, we show the existence of PBTD($n$) for $n \in \{9,11,15\}$.
%It is not difficult to prove that there is no PBTD($n$) for $n \le 4$.
%As a consequence, the existence of PBTD($n$) has been completely determined.

%\begin{df}{\rm 
A {\it partitioned balanced tournament design of side $n$}, \textrm{PBTD($n$)}, defined on a $2n$-set $V$,
is an arrangement of the ${2n \choose n}$
distinct unordered pairs of the elements of $V$ into an $n\times (2n-1)$ arrays such that
\begin{enumerate}
\item every element of $V$ is contained in precisely one cell of each column,
\item every element of $V$ is contained in at most two cells of any row,
\item each row contains all $2n$ elements of $V$ in the first $n$ columns, and
\item each row contains all $2n$ elements of $V$ in the last $n$ columns,
\end{enumerate}
%}\end{df}

\noindent
see \cite{hb}.

E. R. Lamken prove the following theorem.

\begin{thm}[\cite{pbtd}]\label{th:Lamken}
There exists a PBTD{\normalfont ($n$)} for $n$ a positive integer, $n \ge 5$, except possibly for $n \in \{9,11,15\}$.
\end{thm}

Let $V$ be a $2n$-set.
A {\it Howell design} $H(s,2n)$ is an $s \times s$ array, $H$, that satisfies the following three conditions
\begin{enumerate}
\item every cell of $H$ is empty or contains an unordered pair of elements from $V$,
\item each element of $V$ occurs in each row and columns of $H$, and
\item each unordered pair of elements from $V$ occurs in at most one cell of $H$,
\end{enumerate}

\noindent
see \cite{howell}.
%A row $r$ is called a {\it factor} if cells in the row $r$ contains all $2n$ elements of $V$.
%In each row of a Howell design $H(s,2n)$, there exist $s$ cells, {\it a factor}, which contains all $2n$ elements of $V$.
%We can refere Definition and the existence results of Howell designs has been completely determined \cite{howell}.
%A Howell design $H(s,2n)$ has the property that in each row there exist $s$ cells, {\it a factor}, which contains all $2n$ elements of $V$.
%For a partitioned balanced tournament design $T$ of side $n$ ,
For $T$ a PBTD($n$), 
let $T^L,T^C$ and $T^R$ be the first $(n-1)$ columns, the $n$-th column and the last $(n-1)$ columns of $T$, respectively. 
Then ($T^L\ T^C$) and  ($T^R\ T^C$) are Howell designs $H(n,2n)$.
These two designs are called {\it almost disjoint}.
Conversely, if there is a pair of almost disjoint Howell designs, then there is a partitioned balanced tournament design.
%If there exist almost disjoint Howell designs, then 
%We use the termonologies about designs in \cite{hb} and \cite{howell}.

By computer calculation, 
we found almost disjoint Howell designs $H(n,2n)$ for $n \in \{9,11,15\}$ in figures~$1,2$ and $3$.
Hence the following theorem holds.
\begin{thm}\label{th:ArayaTokihisa}
Partitioned balanced tournament designs of side $n$ exist for $n \in \{9,11,15\}$.
\end{thm}

%As there is no partitioned balanced tournament design of side $n$ for $n \le 4$, 
It is not difficult to show that there is no PBTD($n$) for $n \le 4$.
Therefore we have the following corollary from Theorem \ref{th:Lamken} and \ref{th:ArayaTokihisa}.
\begin{cor}
%There exists a PBTD{\normalfont ($n$)} for $n$ a positive integer, $n \ge 5$.
There exists a PBTD{\normalfont ($n$)} if and only if $n$ is a positive integer, $n \ge 5$.
\end{cor}

%\begin{landscape}
\begin{figure}
$\begin{array}{|c|c|c|c|c|c|c|c|c|}
\hline
2, 16  &  3, 17  &  4, 6  &  5, 7  &  8, 10  &  9, 11  &  12, 14  &  13, 15  &  0, 1\\\hline
0, 4  &  1, 5  &  7, 9  &  6, 8  &  11, 13  &  10, 12  &  15, 17  &  14, 16  &  2, 3 \\ \hline
1, 3  &  0, 2  &  10, 13  &  11, 12  &  14, 17  &  15, 16  &  6, 9  &  7, 8  &  4, 5\\ \hline
10, 14  &  11, 15  &  0, 8  &  1, 9  &  2, 4  &  3, 5  &  13, 16  &  12, 17  &  6, 7\\ \hline
5, 6  &  4, 7  &  2, 17  &  3, 16  &  12, 15  &  13, 14  &  0, 10  &  1, 11  &  8, 9\\ \hline
8, 12  &  9, 13  &  1, 15  &  0, 14  &  5, 16  &  4, 17  &  3, 7  &  2, 6  &  10, 11\\\hline
9, 15  &  8, 14  &  11, 16  &  10, 17  &  3, 6  &  2, 7  &  1, 4  & 0, 5  &  12, 13\\ \hline
11, 17  &  10, 16  &  5, 12  &  4, 13  &  1, 7  &  0, 6  &  2, 8  &  3, 9  &  14, 15\\ \hline
7, 13  &  6, 12  &  3, 14  &  2, 15  &  0, 9  &  1, 8  &  5, 11  &  4, 10  &  16, 17 \\\hline
\end{array}$

\smallskip
$\begin{array}{|c|c|c|c|c|c|c|c|c|}
\hline
2, 5  &  3, 4  &  6, 15  &  7, 14  &  8, 11  &  9, 10  &  12, 16  &  13, 17 & 0,1\\\hline
0, 16  &  1, 17  &  4, 8  &  5, 9  &  6, 13  &  7, 12  &  10, 15  &  11, 14 & 2,3\\ \hline
6, 10  &  7, 11  &  1, 16  &  0, 17  &  9, 12  &  8, 13  &  2, 14  &  3, 15 & 4,5\\ \hline
3, 13  &  2, 12  &  9, 17  &  8, 16  &  4, 14  &  5, 15  &  0, 11  &  1, 10 & 6,7\\ \hline
4, 11  &  5, 10  &  2, 13  &  3, 12  &  0, 15  &  1, 14  &  7, 17  &  6, 16 & 8,9\\ \hline
1, 12  &  0, 13  &  5, 14  &  4, 15  &  7, 16  &  6, 17  &  3, 8  &  2, 9 & 10,11\\\hline
9, 14  &  8, 15  &  3, 11  &  2, 10  &  5, 17  &  4, 16  &  1, 6  &  0, 7 & 12,13\\ \hline
8, 17  &  9, 16  &  7, 10  &  6, 11  &  1, 2  &  0, 3  &  5, 13  &  4, 12 & 14,15\\ \hline
7, 15  &  6, 14  &  0, 12  &  1, 13  &  3, 10  &  2, 11  &  4, 9  &  5, 8 & 16,17\\\hline
\end{array}$

\caption{a pair of almost disjoint Howells designs $H(9,18)$}
\end{figure}
%\end{landscape}

\begin{figure}
$\begin{array}{|c|c|c|c|c|c|c|c|c|c|c|}
\hline
 2, 4 & 3, 5 & 18, 21 & 19, 20 & 15, 17 & 14, 16 & 11, 13 & 10, 12 & 7, 8 & 6, 9 & 0, 1 \\\hline 
 0, 9 & 1, 8 & 4, 6 & 5, 7 & 13, 20 & 12, 21 & 17, 19 & 16, 18 & 11, 15 & 10, 14 & 2, 3 \\\hline 
 11, 17 & 10, 16 & 1, 2 & 0, 3 & 6, 8 & 7, 9 & 12, 15 & 13, 14 & 19, 21 & 18, 20 & 4, 5 \\\hline 
 13, 21 & 12, 20 & 11, 19 & 10, 18 & 3, 4 & 2, 5 & 0, 8 & 1, 9 & 14, 17 & 15, 16 & 6, 7 \\\hline 
 16, 19 & 17, 18 & 13, 15 & 12, 14 & 11, 21 & 10, 20 & 5, 6 & 4, 7 & 0, 2 & 1, 3 & 8, 9 \\\hline 
 5, 12 & 4, 13 & 7, 14 & 6, 15 & 9, 16 & 8, 17 & 1, 18 & 0, 19 & 3, 20 & 2, 21 & 10, 11 \\\hline 
 8, 10 & 9, 11 & 17, 20 & 16, 21 & 14, 18 & 15, 19 & 2, 7 & 3, 6 & 1, 5 & 0, 4 & 12, 13 \\\hline 
 3, 7 & 2, 6 & 0, 10 & 1, 11 & 12, 19 & 13, 18 & 16, 20 & 17, 21 & 4, 9 & 5, 8 & 14, 15 \\\hline 
 1, 6 & 0, 7 & 5, 9 & 4, 8 & 2, 10 & 3, 11 & 14, 21 & 15, 20 & 12, 18 & 13, 19 & 16, 17 \\\hline 
 14, 20 & 15, 21 & 3, 8 & 2, 9 & 1, 7 & 0, 6 & 4, 10 & 5, 11 & 13, 16 & 12, 17 & 18, 19 \\\hline 
 15, 18 & 14, 19 & 12, 16 & 13, 17 & 0, 5 & 1, 4 & 3, 9 & 2, 8 & 6, 10 & 7, 11 & 20, 21 \\\hline
\end{array}$

\smallskip
$\begin{array}{|c|c|c|c|c|c|c|c|c|c|c|}
\hline
 2, 12 & 11, 18 & 8, 14 & 4, 15 & 3, 13 & 5, 21 & 9, 17 & 6, 20 & 10, 19 & 7, 16 & 0, 1 \\\hline 
 5, 15 & 4, 14 & 11, 20 & 0, 16 & 6, 17 & 9, 18 & 7, 13 & 1, 19 & 8, 12 & 10, 21 & 2, 3 \\\hline 
 8, 19 & 7, 17 & 6, 16 & 11, 12 & 2, 18 & 10, 13 & 1, 20 & 9, 15 & 3, 21 & 0, 14 & 4, 5 \\\hline 
 4, 20 & 0, 21 & 9, 19 & 8, 18 & 11, 14 & 2, 16 & 10, 15 & 3, 12 & 1, 17 & 5, 13 & 6, 7 \\\hline 
 11, 16 & 6, 12 & 2, 13 & 1, 21 & 0, 20 & 7, 15 & 4, 18 & 10, 17 & 5, 14 & 3, 19 & 8, 9 \\\hline 
 7, 18 & 9, 20 & 1, 12 & 3, 14 & 5, 16 & 6, 19 & 8, 21 & 0, 13 & 2, 15 & 4, 17 & 10, 11 \\\hline 
 9, 14 & 3, 15 & 4, 21 & 2, 17 & 1, 10 & 0, 11 & 5, 19 & 8, 16 & 7, 20 & 6, 18 & 12, 13 \\\hline 
 3, 10 & 1, 16 & 5, 17 & 6, 13 & 4, 19 & 8, 20 & 2, 11 & 7, 21 & 0, 18 & 9, 12 & 14, 15 \\\hline 
 6, 21 & 5, 10 & 3, 18 & 7, 19 & 8, 15 & 1, 14 & 0, 12 & 4, 11 & 9, 13 & 2, 20 & 16, 17 \\\hline 
 0, 17 & 8, 13 & 7, 10 & 5, 20 & 9, 21 & 4, 12 & 3, 16 & 2, 14 & 6, 11 & 1, 15 & 18, 19 \\\hline 
 1, 13 & 2, 19 & 0, 15 & 9, 10 & 7, 12 & 3, 17 & 6, 14 & 5, 18 & 4, 16 & 8, 11 & 20, 21 \\\hline
\end{array}$
\caption{a pair of almost disjoint Howell designs $H(11,22)$}
\end{figure}

\begin{landscape}
\begin{figure}
$\begin{array}{|c|c|c|c|c|c|c|c|c|c|c|c|c|c|c|}
\hline
 2, 4  &  3, 5  &  6, 9  &  7, 8  &  10, 14  &  11, 15  &  12, 16  &  13, 17  &  18, 20  &  19, 21  &  22, 26  &  23, 27  &  24, 29  &  25, 28  &  0, 1 \\\hline 
 9, 10  &  8, 11  &  12, 14  &  13, 15  &  0, 18  &  1, 19  &  24, 28  &  25, 29  &  17, 26  &  16, 27  &  21, 23  &  20, 22  &  4, 6  &  5, 7  &  2, 3 \\\hline 
 1, 15  &  0, 14  &  2, 20  &  3, 21  &  16, 26  &  17, 27  &  23, 25  &  22, 24  &  6, 8  &  7, 9  &  18, 29  &  19, 28  &  11, 12  &  10, 13  &  4, 5 \\\hline 
 5, 23  &  4, 22  &  18, 28  &  19, 29  &  25, 27  &  24, 26  &  17, 20  &  16, 21  &  0, 13  &  1, 12  &  9, 11  &  8, 10  &  3, 15  &  2, 14  &  6, 7 \\\hline 
 17, 21  &  16, 20  &  27, 29  &  26, 28  &  19, 22  &  18, 23  &  11, 13  &  10, 12  &  5, 15  &  4, 14  &  0, 3  &  1, 2  &  7, 25  &  6, 24  &  8, 9 \\\hline 
 16, 28  &  17, 29  &  21, 24  &  20, 25  &  1, 13  &  0, 12  &  2, 5  &  3, 4  &  9, 27  &  8, 26  &  6, 14  &  7, 15  &  19, 23  &  18, 22  &  10, 11 \\\hline 
 22, 27  &  23, 26  &  1, 3  &  0, 2  &  4, 7  &  5, 6  &  8, 14  &  9, 15  &  21, 25  &  20, 24  &  10, 28  &  11, 29  &  16, 18  &  17, 19  &  12, 13 \\\hline 
 0, 19  &  1, 18  &  13, 16  &  12, 17  &  11, 28  &  10, 29  &  9, 26  &  8, 27  &  4, 23  &  5, 22  &  7, 24  &  6, 25  &  2, 21  &  3, 20  &  14, 15 \\\hline 
 13, 20  &  12, 21  &  25, 26  &  24, 27  &  3, 6  &  2, 7  &  10, 18  &  11, 19  &  22, 28  &  23, 29  &  1, 5  &  0, 4  &  9, 14  &  8, 15  &  16, 17 \\\hline 
 26, 29  &  27, 28  &  5, 8  &  4, 9  &  12, 20  &  13, 21  &  3, 7  &  2, 6  &  11, 14  &  10, 15  &  17, 25  &  16, 24  &  1, 22  &  0, 23  &  18, 19 \\\hline 
 6, 11  &  7, 10  &  0, 22  &  1, 23  &  5, 9  &  4, 8  &  19, 27  &  18, 26  &  3, 24  &  2, 25  &  12, 15  &  13, 14  &  17, 28  &  16, 29  &  20, 21 \\\hline 
 3, 25  &  2, 24  &  7, 11  &  6, 10  &  21, 29  &  20, 28  &  0, 15  &  1, 14  &  16, 19  &  17, 18  &  4, 27  &  5, 26  &  8, 13  &  9, 12  &  22, 23 \\\hline 
 8, 12  &  9, 13  &  17, 23  &  16, 22  &  2, 15  &  3, 14  &  6, 29  &  7, 28  &  1, 10  &  0, 11  &  19, 20  &  18, 21  &  5, 27  &  4, 26  &  24, 25 \\\hline 
 18, 24  &  19, 25  &  4, 15  &  5, 14  &  8, 17  &  9, 16  &  21, 22  &  20, 23  &  7, 29  &  6, 28  &  2, 13  &  3, 12  &  0, 10  &  1, 11  &  26, 27 \\\hline 
 7, 14  &  6, 15  &  10, 19  &  11, 18  &  23, 24  &  22, 25  &  1, 4  &  0, 5  &  2, 12  &  3, 13  &  8, 16  &  9, 17  &  20, 26  &  21, 27  &  28, 29 \\\hline
\end{array}$

\smallskip
$\begin{array}{|c|c|c|c|c|c|c|c|c|c|c|c|c|c|c|}
\hline
   2, 11 & 3, 10 & 9, 24 & 15, 27 & 4, 18 & 5, 19 & 16, 23 & 12, 29 & 14, 26 & 7, 20 & 13, 25 & 8, 28 & 17, 22 & 6, 21 & 0, 1 \\\hline
   18, 25 & 4, 13 & 5, 12 & 11, 26 & 15, 29 & 6, 20 & 7, 21 & 8, 23 & 0, 17 & 14, 28 & 9, 22 & 1, 27 & 10, 16 & 19, 24 & 2, 3 \\\hline
   9, 23 & 20, 27 & 1, 6 & 0, 7 & 13, 28 & 15, 17 & 8, 22 & 21, 26 & 10, 25 & 2, 19 & 14, 16 & 11, 24 & 3, 29 & 12, 18 & 4, 5 \\\hline
   10, 24 & 11, 25 & 22, 29 & 3, 8 & 2, 9 & 1, 16 & 15, 19 & 0, 20 & 23, 28 & 12, 27 & 4, 21 & 14, 18 & 13, 26 & 5, 17 & 6, 7 \\\hline
   15, 21 & 12, 26 & 13, 27 & 17, 24 & 5, 10 & 4, 11 & 3, 18 & 7, 19 & 2, 22 & 16, 25 & 0, 29 & 6, 23 & 14, 20 & 1, 28 & 8, 9 \\\hline
   5, 20 & 15, 23 & 0, 28 & 1, 29 & 19, 26 & 7, 12 & 6, 13 & 3, 16 & 9, 21 & 4, 24 & 18, 27 & 2, 17 & 8, 25 & 14, 22 & 10, 11 \\\hline
   1, 8 & 7, 22 & 15, 25 & 2, 16 & 3, 17 & 21, 28 & 0, 9 & 14, 24 & 5, 18 & 11, 23 & 6, 26 & 20, 29 & 4, 19 & 10, 27 & 12, 13 \\\hline
   0, 26 & 2, 28 & 4, 16 & 6, 18 & 8, 20 & 10, 22 & 12, 24 & 1, 25 & 3, 27 & 5, 29 & 7, 17 & 9, 19 & 11, 21 & 13, 23 & 14, 15 \\\hline
   4, 29 & 1, 21 & 7, 18 & 13, 22 & 14, 25 & 2, 23 & 10, 20 & 11, 27 & 6, 19 & 3, 9 & 12, 28 & 15, 26 & 5, 24 & 0, 8 & 16, 17 \\\hline
   12, 22 & 6, 17 & 3, 23 & 9, 20 & 1, 24 & 14, 27 & 4, 25 & 2, 10 & 13, 29 & 8, 21 & 5, 11 & 0, 16 & 15, 28 & 7, 26 & 18, 19 \\\hline
   6, 27 & 0, 24 & 8, 19 & 5, 25 & 11, 22 & 3, 26 & 14, 29 & 9, 28 & 4, 12 & 1, 17 & 10, 23 & 7, 13 & 2, 18 & 15, 16 & 20, 21 \\\hline
   14, 17 & 8, 29 & 2, 26 & 10, 21 & 7, 27 & 13, 24 & 5, 28 & 15, 18 & 11, 16 & 0, 6 & 3, 19 & 12, 25 & 1, 9 & 4, 20 & 22, 23 \\\hline
   7, 16 & 14, 19 & 10, 17 & 4, 28 & 12, 23 & 9, 29 & 1, 26 & 6, 22 & 15, 20 & 13, 18 & 2, 8 & 5, 21 & 0, 27 & 3, 11 & 24, 25 \\\hline
   3, 28 & 9, 18 & 14, 21 & 12, 19 & 6, 16 & 0, 25 & 11, 17 & 5, 13 & 8, 24 & 15, 22 & 1, 20 & 4, 10 & 7, 23 & 2, 29 & 26, 27 \\\hline
   13, 19 & 5, 16 & 11, 20 & 14, 23 & 0, 21 & 8, 18 & 2, 27 & 4, 17 & 1, 7 & 10, 26 & 15, 24 & 3, 22 & 6, 12 & 9, 25 & 28, 29 \\\hline
\end{array}$
\caption{a pair of almost disjoint Howell designs $H(15,30)$}
\end{figure}
\end{landscape}

%\section{Graph representations of a PBTD($n$)}
\section{Observations}
Let $V=\{0,1,\dots,2n-1\}$ be a $2n$-set and $T=(T^L\ T^C\ T^R)$ a PBTD($n$).
Suppose $A$ is the array obtained by permuting elements of $V$, the rows, the first $n-1$ columns, 
the last $n-1$ columns of $T$, or $A=(T^R\ T^C\ T^L)$.
Then $A$ is also a PBTD($n$).
Two PBTD($n$) are {\it isomorphic} if one can be obtained from the other by these operations.
By permuting elements of $V$, we may assume $T^C$ is the transposed of the array
$(\{0,1\}\ \{2,3\}\ \dots\ \{2n-2,2n-1\})$.

From Dinitz and Dinitz~\cite{pbtd10}, there exist two PBTD($5$)'s up to isomorphism.
For these two PBTD($5$)'s, we find that there exists the permutation 
$$\sigma=(0,1)(2,3)(4,5)(6,7)(8,9)$$ such that

$$T^L=
\begin{array}{|c|c|c|c|}
\hline
t_{11} & \sigma(t_{11}) & t_{13} & \sigma(t_{13})\\\hline
t_{21} & \sigma(t_{21}) & t_{23} & \sigma(t_{23})\\\hline
t_{31} & \sigma(t_{31}) & t_{33} & \sigma(t_{33})\\\hline
t_{41} & \sigma(t_{41}) & t_{43} & \sigma(t_{43})\\\hline
t_{51} & \sigma(t_{51}) & t_{53} & \sigma(t_{53})\\\hline
\end{array}
\text{  and  } T^R=
\begin{array}{|c|c|c|c|}
\hline
t_{16} & \sigma(t_{16}) & t_{18} & \sigma(t_{18})\\\hline
t_{26} & \sigma(t_{26}) & t_{28} & \sigma(t_{28})\\\hline
t_{36} & \sigma(t_{36}) & t_{38} & \sigma(t_{38})\\\hline
t_{46} & \sigma(t_{46}) & t_{48} & \sigma(t_{48})\\\hline
t_{56} & \sigma(t_{56}) & t_{58} & \sigma(t_{58})\\\hline
\end{array}
\ .$$
Thus we observe that these two PBTD($5$)'s are determined by some $4$ columns and the permutation $\sigma$.

Seah and Stinson~\cite{pbtd7} obtained two almost disjoint Howell designs $H(7,14)$ 
by computer calcuation for the given $T^L$
which was constructed by E. R. Lamken.
Then for these two PBTD($7$)'s,
we find that there exists the permutation $$\sigma=(0,1)(2,3)(4,5)(6,7)(8,9)(10,11)(12,13)$$ such that 
$$T^L=
\begin{array}{|c|c|c|c|c|c|}
\hline
t_{11} & \sigma(t_{11}) & t_{13} & \sigma(t_{13}) & t_{15} & \sigma(t_{15})\\\hline
t_{21} & \sigma(t_{21}) & t_{23} & \sigma(t_{23}) & t_{25} & \sigma(t_{25})\\\hline
t_{31} & \sigma(t_{31}) & t_{33} & \sigma(t_{33}) & t_{35} & \sigma(t_{35})\\\hline
t_{41} & \sigma(t_{41}) & t_{43} & \sigma(t_{43}) & t_{45} & \sigma(t_{45})\\\hline
t_{51} & \sigma(t_{51}) & t_{53} & \sigma(t_{53}) & t_{55} & \sigma(t_{55})\\\hline
t_{61} & \sigma(t_{61}) & t_{63} & \sigma(t_{63}) & t_{65} & \sigma(t_{65})\\\hline
t_{71} & \sigma(t_{71}) & t_{73} & \sigma(t_{73}) & t_{75} & \sigma(t_{75})\\\hline
\end{array}\ .$$

\noindent
Also we find that there exists the permutation 
$$\tau = (0,2,4)(1,3,5)(8,10,12)(9,11,13)$$ such that

$$T^L=
\begin{array}{|c|c|c|c|c|c|}
\hline
t_{11} & t_{12} & t_{13} & t_{14} & t_{15} & t_{16}\\\hline
\tau(t_{15}) & \tau(t_{16}) & \tau(t_{11}) & \tau(t_{12}) & \tau(t_{13}) & \tau(t_{14})\\\hline
\tau^2(t_{13}) & \tau^2(t_{14}) &\tau^2(t_{15}) & \tau^2(t_{16}) & \tau^2(t_{11}) & \tau^2(t_{12}) \\\hline
t_{41} & t_{42} & \tau(t_{41}) & \tau(t_{42}) & \tau^2(t_{41}) & \tau^2(t_{42})\\\hline
t_{51} & t_{52} & t_{53} & t_{54} & t_{55} & t_{56}\\\hline
\tau(t_{55}) & \tau(t_{56}) & \tau(t_{51}) & \tau(t_{52}) & \tau(t_{53}) & \tau(t_{54})\\\hline
\tau^2(t_{53}) & \tau^2(t_{54}) &\tau^2(t_{55}) & \tau^2(t_{56}) & \tau^2(t_{51}) & \tau^2(t_{52})\\\hline
\end{array}$$
and 
$$T^R=
\begin{array}{|c|c|c|c|c|c|}
\hline
t_{1,8} & t_{1,9} & t_{1,10} & t_{1,11} & t_{1,12} & t_{1,13}\\\hline
\tau(t_{1,10}) & \tau(t_{1,8}) & \tau(t_{1,9}) & \tau(t_{1,13}) & \tau(t_{1,11}) & \tau(t_{1,12})\\\hline
\tau^2(t_{1,9}) & \tau^2(t_{1,10}) & \tau^2(t_{1,8}) & \tau^2(t_{1,12}) & \tau^2(t_{1,13}) & \tau^2(t_{1,11})\\\hline
t_{4,8} & \tau(t_{4,8}) & \tau^2(t_{4,8}) & t_{4,11} & \tau^2(t_{4,11}) & \tau^2(t_{4,11})\\\hline
t_{5,8} & t_{5,9} & t_{5,10} & t_{5,11} & t_{5,12} & t_{5,13}\\\hline
\tau(t_{5,10}) & \tau(t_{5,8}) & \tau(t_{5,9}) & \tau(t_{5,13}) & \tau(t_{5,11}) & \tau(t_{5,12})\\\hline
\tau^2(t_{5,9}) & \tau^2(t_{5,10}) & \tau^2(t_{5,8}) & \tau^2(t_{5,12}) & \tau^2(t_{5,13}) & \tau^2(t_{5.11})\\\hline
\end{array}\ .$$
Thus we observe that $T^L$ is determined by some $7$ cells, the permutations $\sigma$ and $\tau$.
And $T^R$ is determined by some $14$ cells and the permutation $\tau$.

From these observations, we make GAP programs, see~\cite{GAP4}, to construct partitioned balanced tournament designs.
%Note that there are PBTD($6$)'s which have above permutations from master thesis of second author.
And we found designs in figures 1,2 and 3.

\
\section*{acknowledgements}
%\begin{acknowledgement}
This work was supported by JSPS KAKENHI Grant Number 21K03350. 
In this research work we used the supercomputer of ACCMS, Kyoto University.
%\end{acknowledgement}

\end{document}